\documentclass[a4paper,10pt]{article}
\usepackage{amsfonts,latexsym,amssymb}

\renewcommand{\le}{\leqslant} \renewcommand{\ge}{\geqslant}

\newtheorem{theorem}{Theorem}

\newtheorem{theoremno}{Theorem}

\newtheorem{lemma}{Lemma}
\newtheorem{definition}[lemma]{Definition}
\newtheorem{proposition}[lemma]{Proposition}

\begin{document}
\title{Maximal generalization of Baum-Katz theorem
and optimality of sequential tests}
\author{Didier Piau\\ \em Universit\'e Lyon 1}
\date{October 4th, 2000}

\thispagestyle{empty}

\maketitle

\begin{abstract}
Baum-Katz theorem asserts that the Ces\`aro means of i.i.d.\ increments
distributed like $X$ $r$-converge if and only if $|X|^{r+1}$ is
integrable.  We generalize this, and we unify other results, by
proving that the following equivalence holds, if and only if $G$ is
moderate: the Ces\`aro means $G$-converge if and only if
$G(L_{a})$ is integrable for every $a$ if
and only if $|X|\,G(|X|)$ is integrable.  Here,
$L_{a}$ is the last time when the deviation of the Ces\`aro mean
from its limit exceeds $a$, and $G$-convergence is
the analogue of $r$-convergence.  This solves a question about the
asymptotic optimality of Wald's sequential tests.
\end{abstract}

\bigskip



\section{Introduction}
\label{s:intro}

In this paper, we prove a refinement of the strong law of large
numbers and we apply it to a problem of sequential analysis.  Let $X$
and $(X_n)_{n\ge 1}$ denote i.i.d.\ random variables defined on the same
probability space, and
$$
S_n:=X_1+\cdots+X_n.
$$
The most standard version of the strong law of
large numbers asserts that $S_n/n$ converges almost surely if and only if $X$
is integrable.  This amounts to two
statements: first, when $X$ is integrable, the last time
$L^x_{a}$ of a deviation of size at least $a$
of $S_n/n$ from $x=E[X]$ is almost surely finite, for any $a$; and second,
when $X$ is not integrable, no value of $x$ is such that
$L^x_{a}$ is a.s.\ finite for any $a$.
Strassen~(1967) introduced refined versions of this, see also
Lai~(1976), using the notion of $r$-convergence, defined below.

\begin{definition}[Strassen]\label{def1}
Let $r>0$ and $(Y_n)_n$ be random variables defined on the same
probability space.  Set
$$
L_{a}^x:=\sup\{0\}\cup\{n\ge 1\,;\,\|Y_n-x\|\ge a\}.
$$
The sequence $(Y_n)_n$ is said to $r$-converge to $x$ if and only if
$(L_{a}^x)^r$
is integrable, for every positive $a$.
\end{definition}

Complete convergence to $x$, see Hsu and Robbins (1947) and Erd\"{o}s
(1949), is the integrability of the total time spent outside of any
ball around $x$.  Since this time is at most $L^x_{a}$,
$1$-convergence implies complete convergence.  The following
characterizes the $r$-convergence of Ces\`aro means of i.i.d.\
increments.

\begin{theoremno}[Baum and Katz (1965)]
\label{t.bk}
For each positive $r$, $S_n/n$ $r$-converges if and only if $|X|^{r+1}$ is
integrable.  Then, $S_n/n$ $r$-converges to $E[X]$.
\end{theoremno}

We generalize Baum-Katz theorem to every function $G$ such that an
analogue of the result can hold.  We prove that these are exactly the
moderate functions.  This allows to solve an open question of
sequential analysis.

\begin{definition}[Feller~(1969)]\label{def2}
A positive function $G$ defined on $t\ge 0$, is moderate if $G$
is non decreasing, increases to infinity, and if $G(2t)/G(t)$ is
bounded.
\end{definition}

Bingham, Goldie and Teugels~(1989) call the moderate functions,
``increasing functions of dominated variation'' (see their
section~2.1).  For any deterministic non decreasing function $G$,
increasing to infinity, $(Y_n)_n$ $G$-converges to $x$ if and only
if $G(L_{a}^x)$ is integrable for every $a$.
Here is our main result.

\begin{theorem}\label{th0}
Let $X$ and $(X_n)_{n\ge 1}$ denote i.i.d.\ random variables,
defined on the same
 probability space, and let $G$ denote a moderate function.  Then, $S_n/n$
 $G$-converges if and only if $|X|\,G(|X|)$ is integrable.  More
 precisely, assume that $X$ is integrable and centered. Then,
 assertions {\bf (a)}, {\bf (b)} and {\bf (c)} below are equivalent. 

\hspace*{2ex}
{\bf (a)} $|X|\,G(|X|)$ is integrable.

\hspace*{2ex} 
{\bf (b)} For every $a$,
$S(X,G,a):=\sum_{n\ge
1}n^{-1}\,G(n)\,P[|S_n/n|\ge a]$ is finite.

\hspace*{2ex}
{\bf (c)}
For every $a$,
$G(L_{a}^0)$ is integrable.

Furthermore, if $G$ is non decreasing, unbounded, and 
non moderate,
assertions {\bf (a)} and {\bf (c)} are not equivalent.
\end{theorem}
\begin{theorem}
\label{c.cor}
Let $G$ denote a moderate function.  Under an integrability assumption
analogous to assumption {\bf (a)} in theorem~\ref{th0}, 
Wald's sequential test is asymptotically optimal,
when the error probabilities goes to zero, with respect to the
$G$-moment of the observation time of the test.
\end{theorem}

See section~\ref{s.sa} for a precise statement of theorem~\ref{c.cor}.
Baum and Katz~(1965) prove theorem~\ref{th0} for $G(t)=t^r$ with
$r>0$, as a part of their theorem~3, and they prove a similar statement for
$G(t)=t^r\,\log(t)$ with $r\ge 0$, as a part of their theorem~2.
All these functions $G$ are moderate.  Baum and Katz~(1965) use results of
Katz~(1963). Some methods of these two papers are from
Erd\"{o}s~(1949).  We use some methods of these three papers but none
of their results.

Lanzinger~(1998) studies a formally equivalent extension of Baum-Katz
law. Basically, he shows that, in our theorem~\ref{th0}, {\bf (a)} holds for
$G(t)=\exp(bt)$ if and only if {\bf (b)} holds for
$G(t)=\exp(ct)$, and for $a$ large enough (namely,
for every $a>c/b$).  Thus, Lanzinger's result applies to
a law of large numbers if $\exp(c\,|X|)$ is integrable for every positive
$c$.  This condition is not equivalent to the integrability of
$H(|X|)$, for any single function $H$.  (To see this, note that,
if $H$ is a
candidate, then $t=o(\log H(t))$ when $t\to\infty$. Choose $\Omega$ such
that $\Omega=o(H)$ and such that $t=o(\log\Omega(t))$.  Then, $\Omega(|X|)$
may be integrable while $H(|X|)$ is not.)

The paper is organized as follows.  In section~\ref{s.sa}, we deal with the
application of theorem~\ref{th0} to sequential analysis.  In
section~\ref{s.auxi}, we reduce theorem~\ref{th0} to the case of
symmetric laws, and we state explicit universal bounds of the
quantities involved in assertions {\bf (a)},
{\bf (b)}
and {\bf (c)} in terms of the quantities involved in assertions 
{\bf (c)}, {\bf (a)} and {\bf (b)} 
respectively, for symmetric
laws.  The proofs of these explicit bounds are rather cumbersome, and
we only sketch them, in section~\ref{s.proofs}.  This
section also proves theorem~\ref{th0} for non moderate functions.  An
appendix gives detailed proofs.


\section{Sequential analysis}
\label{s.sa}

Theorem~\ref{c.cor} answers a question of Koell~(1995) about the
asymptotic optimality of Wald's sequential tests.  Koell proves that
{\bf (c)} implies {\bf (a)} for $a=1/2$ when stronger assumptions hold
(his result is a version of proposition~\ref{th1} below).  The main
result of Koell~(1995) generalizes Lai~(1981) and is as follows.

Assume that the unknown law of the process $(Y_n)_n$ is one of the
finitely many distinct probability measures $P_i$, and call $H_i$ the
hypothesis
$$
H_i:=[\,\mbox{The law of\ }(Y_n)_n\ \mbox{is}\ P_i\,].
$$ 
Denote by ${\mathcal F}_n$ the $\sigma$-field generated by
$(Y_k)_{k\in[1,n]}$, choose a probability $Q$ such that each $P_i$ is
absolutely continuous with respect to $Q$, and a probability $P$ which
is locally equivalent to each $P_i$.  Call $p_{in}$, resp. $p_n$, the
density of $P_i$, resp. $P$, with respect to $Q$, when restricted to
${\mathcal F}_n$.  The likelihood ratio of the sample
$(Y_k)_{1\le k\le n}$ under hypothesis $H_i$ is
$$
R^i_n:=\frac{p_n(Y_1,\ldots,Y_n)}{p_{in}(Y_1,\ldots,Y_n)}.
$$
A decision rule is a couple $(\tau,d)$, where the observation time
$\tau$ is an integrable $({\mathcal F}_n)_n$ stopping time and the
decision $d$ is a ${\mathcal F}_{\tau}$ random variable with values in
$\{H_i\}_i$.  Given error probabilities $a=(a_i)_i$, a stopping time
$\tau$ belongs to the class $T(a)$ if there exists a decision rule
$(\tau, d)$ such that $P_i[d\neq H_i]\le a_i$ for every $i$.

Wald's sequential tests are commonly used decision rules, defined as
follows.  Fix a (multi)level $c=(c_i)_i$, and denote by $\varrho_i$ the
rejecting time of $H_i$ at level $c_i$, that is,
$$
\varrho_i:=\inf\{n\ge 1\,;\,R^i_n\ge c_i\}.
$$
Wald's sequential test at level $c$ is defined by the relations
$$\tau_c:=\min_i\max_{j\neq
i}\varrho_j,\quad
d_c:=H_k,\quad k:=\mbox{argmax}_j\varrho_j.
$$  
Theorem~2 of
Koell (1995) proves the asymptotic optimality of $(\tau_c,d_c)$ with
respect to $E[\tau^r]$.  More precisely, assume that, for each $i$, $
n^{-1}\,\log R^i_n $ $r$-converges under $P$. Then, for each
(multi)probability of error $a=(a_i)_i$, there exists a level
$c=(c_i)_i$ such that $(\tau_c,d_c)$ belongs to $T(a)$. Furthermore,
$c$ can be chosen such that $\tau_c$ is asymptotically the smallest
observation time in $T(a)$ when $a\to 0$, in other words,
$$
E_i[\tau_c^r]=\inf\{E_i[\tau^r]\,;\,\tau\in T(a)\}\cdot(1+o(1)),
$$
for all $i$, when $a\to 0$.  Lai (1981) showed this for two laws $P_1$
and $P_2$.  Koell also proves this, if $t\mapsto t^r$ is replaced by
an increasing function $G$ in ${\mathcal A}$, where ${\mathcal A}$ is
a specific strict subset of the space of continuous moderate
functions, defined in Kroell~(1985).  Our
theorem implies that the same result is valid for any moderate
function.

As an example, assume that the function $G$ is moderate, that the process
$(Y_n)_n$ is i.i.d., and that
$$
\log[p_1(Y_1)/p_{i1}(Y_1)]\,G(\log[p_1(Y_1)/p_{i1}(Y_1)])
$$
is integrable for every $i$.  Then, $(\tau_c,d_c)$ is asymptotically
optimal with respect to $E[G(\tau)]$, that is,
$$
E_i[G(\tau_c)]=\inf\{E_i[G(\tau)]\,;\,\tau\in T(a)\}\cdot(1+o(1)),
$$
for every $i$, when $a\to 0$.
The preceding was an incentive to prove theorem~\ref{th0}.

\section{Effective bounds}
\label{s.auxi}
Theorem~\ref{th0} is a consequence of propositions~\ref{th1} to
\ref{th4} below.  Propositions~\ref{th1} to \ref{th3} provide
effective bounds of {\bf (a)}, {\bf (b)} and {\bf (c)} in terms of
{\bf (c)}, {\bf (a)}  and
{\bf (b)} respectively, when the law of $X$ is symmetric.
Proposition~\ref{th4} reduces the general case to the case of
symmetric laws.  We write $L_{a}$ for
$L_{a}^0$.

\begin{proposition}\label{th1}
Let $G$ denote a moderate function with $G(2t)\le c\,G(t)$ for every positive $t$, and
let $X$ denote an integrable random variable.  
For $\alpha\in]0,1[$, let $t$ be such that
$ E[|X|\,;\,|X|\ge t]\le (1-\alpha).  $ Then,
$$
E[|X|\,G(|X|)]\le 4c^2\,\left\{t\,G(t)+
\alpha^{-1}\,E[G(L_{1/2})]\right\}.
$$
\end{proposition}
\begin{proposition}\label{th2}
\addtocounter{equation}{1} Let $p\ge 1$ be an integer, and let $G$ denote a
non decreasing function such that $G(t)/t^{p+1}$ is integrable on
$(1,+\infty)$ (call this the condition~(\theequation)).
Let $H$ be a non decreasing function
such that, for any $n\ge 1$,
\begin{equation}
\label{e.hpsi}
H(n)\ge n^p\,\sum_{k\ge n}G(k)\,k^{-(p+1)}.
\end{equation}
Then, for any symmetric $X$,
$
S(X,G,1)\le
K_p(X)$, where
$$
K_p(X):=E[|X|\,G(|X|)]+
2^{-p}\,(2p)!\,E[1+|X|]^{p-1}\,E[|X|\,H(|X|)].
$$
If $G$ is moderate, (1) holds for $p$ large enough and
$H=c\,G$ satisfies (\ref{e.hpsi}) for $c$ large enough.
\end{proposition}
\begin{proposition}\label{th3}
Let $G$ denote a non decreasing function,
and $X$ denote a symmetric random variable.
Then,
$$
E[G(L_1)]
\le
G(0)+12\,S(X,G,1/8).
$$
\end{proposition}
Denote by $X'$ an independent 
copy of $X$ and call $X^*:=X-{X'}$ a symmetrized version of $X$. 
\begin{proposition}\label{th4}
Let $G$ be moderate.
Each of the assertions {\bf (a)}, {\bf (b)} and {\bf (c)}
of theorem~\ref{th0} holds for $X$
if and only if it holds for the symmetrized $X^*$.
\end{proposition}

From proposition~\ref{th1}, {\bf (c)} implies {\bf (a)} if $G$ is moderate.
From proposition~\ref{th2}, {\bf (a)} implies {\bf (b)} if $G$ is moderate and
if the law of $X$ is symmetric.  To see that the series in {\bf (b)} then
converges for every $a$, note that, if {\bf (a)} holds for $X$,
then {\bf (a)} holds for every $X/a$ as well.  From
proposition~\ref{th3}, {\bf (b)} implies {\bf (c)} if $G$ is non decreasing
and if the law of $X$ is symmetric.  To see that
$G(L_{a})$ is then integrable for any $a$,
note that $L_1$ for the sum of the random variables
$X_n/a$ is $L_{a}$ for $S_n$.

\section{Sketches of proofs}
\label{s.proofs}
The proof of proposition~\ref{th1} uses precise estimates of 
$$P[\sup_{k\ge n}|X_k/k|\ge 1]$$ through 
Poincar\'e's formula for the probability of
a union of sets, 
Abel's transformation for discrete sums, and the moderation of $G$.

In proposition~\ref{th1}, $T_n$ is the sum of $n$
truncated variables $$X_{kn}:=X_k\,\mbox{\bf 1}_{\{|X_k|\le n\}}.$$
A multinomial expansion allows to estimate $E[(T_n)^{2p}]$.
Since $X_{nk}$ is symmetric, the terms with at least one odd power
disappear.
A bound of the remaining terms involves a decomposition along the first
index $i$ such that $|X_{in}|$ is maximal amongst $\{|X_{kn}|\}_{k\le
n}$.
Precise combinatorics then yield the result.

In proposition~\ref{th3}, one decomposes the
sum along the intervals $[2^i,2^{i+1}[$.
From a maximal inequality for symmetric random variables, due to L\'evy,
the events ``$\max|S_n|/n$ is large for $n$ in an interval''  
and ``$|S_n|/n$ is large for the last $n$ in this interval'' have
the same probability, up to a factor of $2$.
One then reverses the time on each interval, uses L\'evy's inequality
again, and another decomposition of a sum along the  intervals
$[2^i,2^{i+1}[$ concludes the proof.

Proposition~\ref{th4} follows from symmetrization techniques and
estimations of the probability of a deviation from any median, 
as exposed in section~18.1 of
Lo\`eve (1977) for example.
Furthermore, if $(X'_n)_n$ is an independent copy of $(X_n)_n$, and if
$(X^*_n)_n$, where $X^*_n:=X_n-X'_n$,  is a symmetrization  of
$(X_n)_n$, one can compare the events $\{L_{a}\ge n\}$, 
$\{L'_{a}\ge n\}$
and
$\{L^*_{a}\ge n\}$, for different values of $a$
and $n$, where $L'_{a}$ and $L'_{a}$ are the
functional $L_{a}$, when applied to $X'$ and $X^*$
instead of $X$.

The proof of theorem~\ref{th0} for non moderate functions is as
follows.  If {\bf (a)} and {\bf (c)} are equivalent, $G(L_{a})$
is integrable for every positive $a$ as soon as $|X|\,G(|X|)$ is
integrable.  But $L_{2a}$ for $(2X_n)_n$ is
$L_{a}$ for $(X_n)_n$, hence $|X|\,G(2|X|)$ should be
integrable as soon as $|X|\,G(|X|)$ is.  
\\ 
Thus, let $G$ be a
non decreasing, increasing to infinity, non moderate function.  Then,
$G(2t_n)\ge n\,G(t_n)$ for an increasing sequence $(t_n)_n$.
Let the law of $X$ have as support the set $\{\pm t_n\,;\,n\ge 1\}$,
and give to $t_n$ and $-t_n$ the same weight of
$c/(n^2\,t_n\,G(t_n))$, for a given constant $c$. The sum of the
weights converges, hence $c$ can be chosen so as to get a probability
measure. Then, $E[|X|\,G(|X|)]$ is finite and $E[|X|\,G(2|X|)]$
is infinite.  Finally, {\bf (a)} and {\bf (c)} are not equivalent for $X$.

\section*{Appendix: Detailed proofs}

\section{Proof of proposition~\ref{th1}}
\label{s:proof1}
Let $U_n:=S_n/n$.
Since $L_a$ is an a.s.\ finite nonnegative integer,
\begin{equation}
\label{e:espl}
E[G(L_{a})] = 
G(0)+\sum_{n\ge 1}(G(n)-G(n-1))\,P[\sup_{k\ge n}|U_k|\ge a].
\end{equation}
Starting from an interval of
length $2a$, for example $]-a,+a[$,
one cannot make a jump whose  length is larger than $2a$ 
and still be in this
interval after the jump. 
Assume from now on that $a=1/2$. The preceding remark gives:
$$
\{\sup_{k\ge n}|U_k|\ge 1/2\}\supset E_n:=\{\sup_{k\ge n}|X_k/k|\ge 1\}.
$$
The idea of the proof is to
estimate precisely $P[E_n]$
as $n\to\infty$.
Denote by
$$
\sigma_k:=\sum_{i\ge k}P[|X_i|\ge i]=\sum_{i\ge k}P[|X|\ge i].
$$
Since $E_n$ is a union of independent events, 
Poincar\'e's formula at the first and second orders gives
$$
\sum_{k\ge n}P[|X_k|\ge k]\,\left(1-\sum_{i>k}P[|X_i|\ge i]\right)
\le
P[E_n] \le\sum_{k\ge n}P[|X_k|\ge k],
$$
hence $\sigma_n\,(1-\sigma_{n+1})\le P[E_n]\le\sigma_n$.
Since
$X$ is integrable,
$\sigma_{n+1}\to 0$ and $P[E_n]\sim\sigma_n$ 
when $n\to\infty$.
More quantitatively, 
assume that $n\ge n_{X}=\lfloor t\rfloor$, where
$t$ is defined in the statement of
proposition~\ref{th1}.
Then, $n+1\ge t$ and
\begin{eqnarray*}
\sigma_{n+1} & \le & E\left[|X|-n\,;\,{|X|}\ge n+1\right]
\\
& \le & E[|X|\,;\,|X|\ge t]\le 1-\alpha.
\end{eqnarray*}
This yields
$
\alpha\,\sigma_n\le P[E_n]
$
for $n\ge n_X$.
When $n<n_{X}$,
$$
\alpha\,\sigma_n\le P[E_n]+\alpha\,(n_{X}-n).
$$
Going back to the equation~(\ref{e:espl}),
a lower bound of $E[G(L_{1/2})]$ is
$$
E[G(L_{1/2})]\ge
G(0)+\sum_{n\ge 1}(G(n)-G(n-1))\,\alpha\,(\sigma_n-(n_{X}-n)^+).
$$
An application of Abel's transformation to the right hand member gives:
\begin{eqnarray}\label{e.ss}
\sum_{n\ge 1}G(n)\,P[|X|\ge
n] & \le & \alpha^{-1}\,E[G(L_{1/2})]+\sum_{n=1}^{n_{X}-1}G(n)
\\
 & \le & \alpha^{-1}\,E[G(L_{1/2})]+(n_{X}-1)\,G(n_{X}-1)
\nonumber
\\
 & \le & \alpha^{-1}\,E[G(L_{1/2})]+(t-1)\,G(t).
\nonumber
\end{eqnarray}
There should be a multiple of $G(0)$ in the left hand side, 
but the coefficient of  $G(0)$ 
is $1-\alpha(\sigma_1-(n_X-1))$, which is nonnegative because:
$$
\sigma_1\le n_X+\sigma_{n_X+1}\le n_X+1-\alpha\le n_X+\alpha^{-1}-1.
$$
Decomposing the sum of (\ref{e.ss})
along the
intervals $I(i):=[2^i,2^{i+1}[$, one gets:
\begin{equation}
\label{e.sum}
\sum_{n\ge 1}G(n)\,P[|X|\ge
n] \ge\sum_{i\ge 0}2^i\,G(2^{i})\,P[|X|\ge 2^{i+1}].
\end{equation}
On the other hand, 
\begin{eqnarray*}
E[|X|\,G(|X|)] & \le & \sum_{n\ge 0}((n+1)\,G(n+1)-n\,G(n))\,P[|X|\ge n]
\\
 & \le &  G(1)+\sum_{i\ge 0}(2^{i+1}G(2^{i+1})-2^iG(2^i))\,P[|X|\ge 2^i].
\end{eqnarray*}
The last sum runs over integers $i\ge 0$.
We replace the term $i=0$ by
its value
 and for any
$i\ge 1$, we bound the coefficient of $P[|X|\ge 2^i]$ by
$2^{i+1}G(2^{i+1})$. 
This yields
\begin{eqnarray*}
E[|X|\,G(|X|)] & \le & 
2\,G(2)+\sum_{i\ge 1}2^{i+1}G(2^{i+1})\,P[|X|\ge 2^i]
\\
 & \le & 2\,G(2)+4c^2 \sum_{i\ge 0}2^iG(2^i)\,P[|X|\ge 2^{i+1}],
\end{eqnarray*}
by an application of the moderation property of $G$ between $2^{i+1}$ and
$2^{i-1}$.
Composing this with the inequalities~(\ref{e.ss}) and (\ref{e.sum}), one gets 
the first assertion of proposition~\ref{th1}, up to an
 additional term in
the right hand side.
This extra term reads as $2\,G(2)-4c^2\,G(t)$, which is
nonpositive and therefore may be omitted.

\section{Proof of proposition~\ref{th2}}
\label{s:proof2}
In this section, the law of $X$ is symmetric.
For $1\le k\le n$, set
$$
X_{kn}:=X_k\,\mbox{\bf 1}_{\{|X_k|\le n\}},\qquad
T_n:=\sum_{k=1}^nX_{kn}.
$$
If  $|X_k|\le n$ for all $1\le k\le n$, then $S_n=T_n$. Hence:
$$
P[|S_n|\ge n]\le n\,P[|X|\ge n]+P[|T_n|\ge n].
$$
We write $P[|X|\ge n]$ as a sum of $P[k\le |X|<k+1]$ and we use the fact that
$G$ is non decreasing. 
This yields
$$
\sum_{n\ge 1}n^{-1}\,G(n)\,n\,P[|X|\ge n]\le E[|X|\,G(|X|)].
$$
For $P[|T_n|\ge n]$,
we use the Markov inequality with an even integer $2p$, 
whose value will be specified
later:
$$
P[|T_n|\ge n]\le n^{-2p}\,E[(T_n)^{2p}].
$$
In the expansion of $(T_n)^{2p}$, all the terms with at least one odd power disappear
because $X$ is symmetric, hence
$$
E[(T_n)^{2p}]=\sum_cM(2p,2c)\,E[X_{1n}^{2c_1}\,X_{2n}^{2c_2}\cdots
X_{nn}^{2c_n}],
$$
where $c=(c_i)_{1\le i\le n}$ runs over all the  nonnegative integer valued
sequences of length $n$ and of
 sum 
$p$, and where $M(2p,2c)$ is the multinomial coefficient
associated to the sequence $2c$. 
Recall that for any $p\ge 0$ and any nonnegative integer valued
sequence $c=(c_i)_{1\le i\le n}$ of sum
 $c_1+c_2+\cdots+c_n=p$, one has
$$
M(p,c):=\frac{p!}{c_1!\,c_2!\,\ldots\,c_n!}.
$$
If $0\le q\le p$, we write $M(p,q):=M(p,(q,p-q))$ for the usual binomial
coefficient.

The value of the expectation $E[X_{1n}^{2c_1}\,X_{2n}^{2c_2}\cdots
X_{nn}^{2c_n}]$ is entirely determined by the sequence $d=(d_i)_{i}$ 
which is made by listing 
the non zero terms of $c$ and by writing them in the same order.
This value is:
\begin{equation}\label{e:end}
e(n,d):=E[X_1^{2d_1}\,X_2^{2d_2}\ldots X_q^{2d_q};\,|X_0|\le n], 
\end{equation}
where we introduced
$$
|X_0|:=\max\{|X_i|\,;\,1\le i\le q\}.
$$
In equation~(\ref{e:end}), the sum of $d$ is $p$, the length of $d$ is $q$,
and all the terms of $d$ are non zero. 
Exactly $M(n,q)$ sequences $c$ give the same sequence $d$ of length
$q$.
Since $M(n,q)\le n^q\,(q!)^{-1}$, this yields
$$
E[(T_n)^{2p}]\le\sum_dM(2p,2d)\,n^q\,(q!)^{-1}\,e(n,d),
$$
where $q$ is the length of the sequence $d$.
Summing up over $n$, one gets 
$$
\sum_{n\ge 1}n^{-1}\,G(n)\,P[|T_n|\ge n]
\le
\sum_dM(2p,2d)\,(q!)^{-1}\,e(d),
$$
where $e(d)$ is the sum
$$
e(d):=\sum_{n\ge 1}G(n)\,n^{q-2p-1}\,e(n,d).
$$
We now bound each $e(d)$.
Summing up the equations~(\ref{e:end}) yields
$$
e(d)=E\left[\prod_{i=1}^q|X_i|^{2d_i}\cdot
\sum_{n\ge |X_0|}G(n)\,n^{q-2p-1}\right].
$$
Let $p$ denote an integer such that the condition of 
proposition~\ref{th2} holds. 
Since $q\le p$ for every $d$, one has
$$
\sum_{n\ge |X_0|}G(n)\,n^{q-2p-1}\le |X_0|^{q-p}\,\sum_{n\ge
|X_0|}
G(n)\,n^{-(p+1)}\le |X_0|^{q-2p}\,
H(|X_0|).
$$
Since $|X_i|^{2d_i-1}\le|X_0|^{2d_i-1}$ for every $1\le i\le q$, one has
$$
e(d)\le
E\left[H(|X_0|)\,\prod_{i=1}^q|X_i|\right].
$$
Write $i_0$ for anyone of the integers such
that $|X_{i_0}|=|X_0|$,
for example the smallest one.
Then, the expectation of the right hand side, when restricted to the event
$\{i_0=j\}$, is at most
$$
E\left[H(|X_j|)\,|X_j|\,\prod_{i\neq  j}|X_i|\,;\,i_0=j\right]
\le E[|X|]^{q-1}\,E[|X|\,H(|X|)].
$$
There are at most $q$ possible values of $i_0$ and, for a fixed length $q$,
there are $M(p-1,q-1)$ possible values of the $q$-uplet $d$.
(For completeness, we prove this last assertion in 
lemma~\ref{lmpq} at the end
of this section.)
For a given $q$, $M(2p,2d)$ is maximal when each $d_i$ is 
as close as possible of
$p/q$. For $1\le q\le p$, this upper bound is maximal
for $q=p$. Hence, 
$$
M(2p,2d)\le c(p):=(2p)!/2^p.
$$
This shows that
proposition~\ref{th2} holds with the bound $E[|X|\,G(|X|)]+c_p(X)$, where
\begin{eqnarray*}
c_p(X) & := & c(p)
        \sum_{q=1}^pM(p-1,q-1)\,q\,(q!)^{-1}\,E[|X|]^{q-1}\,E[|X|\,H(|X|)]
\\
& \le & c(p)\,E[1+|X|]^{p-1}\,E[|X|\,H(|X|)].
\end{eqnarray*}
The last assertions of proposition~\ref{th2} about moderate functions
are direct consequences of the inequality $G(2t)\le c\,G(t)$.
The combinatorial lemma used during the proof is the following.

\begin{lemma}\label{lmpq}
For $1\le q\le p$,
there are $M(p-1,q-1)$ sequences $d=(d_i)_{1\le i\le q}$, integer valued, 
of length $q$ and sum $p$, such that $d_i\ge 1$ for every $i$.
\end{lemma}
\textbf{Proof of lemma~\ref{lmpq}}
We count the integer valued sequences $b=(b_i)_{1\le i\le q}$ of sum
$(p-q)$ such that $b_i\ge 0$ for every $1\le i\le q$.
The cardinal of this set is the coefficient of $z^{p-q}$ in the expansion:
$$
\prod_{i=1}^q\left(\sum_{n=0}^{+\infty}z^n\right)=(1-z)^{-q},
$$
and this coefficient is $M(p-1,q-1)$.
This ends the proof of lemma~\ref{lmpq}.

\section{Proof of proposition~\ref{th3}}
\label{s:proof3}
In this section, the law of $X$ is symmetric.
In Section~\ref{s:proof1}, (\ref{e:espl}) provides $E[G(L_1)]$
as the sum of a series.
The decomposition of this series
along the intervals $I(i):=[2^i,2^{i+1}[$ yields
\begin{eqnarray*}
E[G(L_{1})] & \le & 
G(0)+\sum_{i\ge 0}P[L_1\ge
2^i]\,\sum_{n\in I(i)}(G(n)-G(n-1))
\\
 & \le & G(0)\,P[L_1=0]+
\sum_{i\ge 0}(G(2^{i+1}-1)-G(2^i-1))\,P[L_1\ge
2^i].
\end{eqnarray*}
The probabilities written in the right hand side are bounded by
$$
P[L_1\ge
2^i]=P[\sup_{n\ge 2^i}|U_n|\ge 1]\le\sum_{j\ge
i}P[\max_{n\in I(j)}|U_n|\ge 1],
$$
where $U_n:=S_n/n$.
Furthermore,
$$
\{\max_{n\in I(j)}|U_n|\ge1\}\subset D_j:=\{\max_{n\in I(j)}|S_n|\ge 2^i\}
$$
An interversion of the order of the summations on $i$ and $j$ yields
\begin{equation}
\label{e.fl1}
E[G(L_1)]\le G(0)+\sum_{i\ge 0}G(2^{i+1})\,P[D_i].
\end{equation}
Since $X$ is symmetric, zero is a median of every $S_n$ and a maximal
inequality due to Paul L\'evy, see Lo\`eve~(1977), section~18.1, states that
$$
P[\max_{n\le m}|S_n|\ge t]\le 2\,P[|S_m|\ge t]
$$
for any $m\ge 1$ and $t\ge 0$.
This yields an upper bound of $P[D_i]$ as follows.
First, $D_i$ can only occur if $|S_{2^{i+1}}|\ge
2^i/2$, or if  $D'_i$ is realized, where
$$
D'_i:=\{\exists n\in I(i)\,;\,|S_n-S_{2^{i+1}}|\ge 2^i/2\}.
$$
We now bound the probability $P[D'_i]$.
The process 
$$
(S_{2^{i+1}}-S_n)_{n\in I(i)}
$$
follows the law of the time reversal of $(S_n)_{1\le n\le 2^i}$, 
hence $P[D'_i]=P[D''_i]$ with
$$
D''_i:=\{\exists n\le 2^i,\,|S_n|\ge 2^i/2\}.
$$
L\'evy's inequality yields
$$
P[D''_i]\le P[\max_{n\le 2^{i+1}}|S_n|\ge 2^i/2]
\le 2\,P[|S_{2^{i+1}}|\ge 2^i/2].
$$
Coming back to the event $D_i$, we proved that
\begin{equation}
\label{e.fl2}
P[D_i]\le 3\,P[|S_{2^{i+1}}|\ge 2^i/2].
\end{equation}
This yields
\begin{eqnarray*}
E[G(L_1)] & \le & 
G(0)+\sum_{i\ge 0}3\,G(2^{i+1})\,P[|U_{2^{i+1}}|\ge1/4]
\\
& = & 
G(0)+\sum_{i\ge 1}3\,G(2^{i})\,P[|U_{2^{i}}|\ge1/4].
\end{eqnarray*}
On the other hand,
we have to estimate 
$$
s:=\sum_{n\ge 1}n^{-1}\,G(n)\,P[|U_n|\ge1/8].
$$
A decomposition of $s$ along the powers of $2$ gives
$$
s\ge\sum_{i\ge 0}2^{-(i+1)}\,G(2^i)\sum_{n\in I(i)}P[|U_n|\ge1/8].
$$
Recall that $X$ is symmetric and let $n$ belong to $I(i)$. Since the
law of  $S_k$  is symmetric, 
\begin{eqnarray*}
P[|S_{2^i}|\ge 2^{i-2}] & \le  & P[\max_{j\le n}|S_j|\ge 2^{i-2}]
\\
 & \le  & 2\,P[|S_n|\ge
2^{i-2}]\le 2\,P[|S_n|\ge n/8].
\end{eqnarray*}
Finally,
\begin{equation}
\label{e.fl3}
4\,s\ge\sum_{i\ge 0}G(2^i)\,P[|U_{2^i}|\ge1/4].
\end{equation}
The inequalities (\ref{e.fl1}), (\ref{e.fl2}) and (\ref{e.fl3}) 
prove the assertion of proposition~\ref{th3}.

\section{Proof of proposition~\ref{th4}}
\label{s:proof4}
We recall some basic facts about the symmetrization
of a random variable around one of its medians.
A median of a real random variable $Y$ is any real number $\mu(Y)$ such that
$$
P[Y\le\mu(Y)]\ge 1/2
\quad
\mbox{and}
\quad
P[Y\ge\mu(Y)]\ge 1/2.
$$
Let $Y'$ be an independent copy of $Y$, independent of all the
other random variables.
Then, the law of $Y^*:=Y-{Y'}$ is symmetric and
\begin{equation}\label{e:sym}
P[|Y-\mu(Y)|\ge 2t]\le 2\,P[|Y^*|\ge 2t]\le 4\,P[|Y|\ge t]
\end{equation}
for every $t$, see Lo\`eve~(1977), section~18.1.
If $H$ is non decreasing, equation~(\ref{e:sym}) yields 
\begin{equation}\label{e:symint}
E[H(|Y-\mu(Y)|)]\le 2\,E[H(|Y^*|)]\le 4\,E[H(2\,|Y|)].
\end{equation}
In order to apply this to $X_n$, we introduce 
$$
\sum_{k=1}^nX^*_k=S_n-S_n',
$$
which is a symmetrized version $S^*_n$ of $S_n$. 
Set $U^*_n:=S^*_n/n$ and let 
$G$ be a moderate function. Then
$$
E[|X^*|\,G(|X^*|)]\le 4\,E[|X|\,G(2|X|)]\le 4\,c\,E[|X|\,G(|X|)].
$$
Hence, (i) for $X$ implies (i) for $X^*$.
For the reverse implication, (\ref{e:symint}) yields
$$
E[|X-m|\,G(|X-m|)]\le 2\,E[|X^*|\,G(|X^*|)],
$$
where $m:=\mu(X)$ is a median of $X$. 
The  expectation of the left hand side
is finite, when it is restricted to the set
$\{|X-m|\ge m\}$, and, on this set, $|X-m|\ge |X|/2$.
Hence,
$$
E[|X|\,G(|X|/2)\,;\,|X-m|\ge m]
$$
is finite.
The complete expectation $E[|X|\,G(|X|/2)]$ is finite as well. 
\\
Finally, 
$G(t)\le c\,G(t/2)$ for any $t$, hence
{\bf (a)} for $X^*$ implies {\bf (a)} for $X$.

The equivalence of {\bf (b)} for $X$ and {\bf (b)} for $X^*$ stems from the
fact that, if $\mu(\xi_n)$ is a median of $\xi_n$ and if $\xi_n$
converges to zero in probability, then $\mu(\xi_n)$ converges to zero.
(To see this, assume that $\xi_n$
converges to zero in probability. Then, 
for every positive $\alpha$, $P[|\xi_n|\ge \alpha]<\frac12$ for $n$
large enough.  
This means that
$|\mu(\xi_n)|\le \alpha$ for $n$ large enough.
Since $\alpha$ is arbitrary, $\mu(\xi_n)$ converges to 
zero.)

Since $X$ is integrable and centered, $U_n$ converges almost surely to zero
and $\mu(U_n)$ converges to zero.
For a fixed positive $a$ and
for $n$ large enough, 
\begin{eqnarray*}
P[|U_n|\ge a] & \le & P[|U_n-\mu(U_n)|\ge a/2]
\\
& \le & 2\,P[|U^*_n|\ge a/2]\le 4\,P[|U_n|\ge a/4],
\end{eqnarray*}
where the convergence of $\mu(U_n)$ to zero yields
the first inequality and (\ref{e:sym}) at the beginning of
this section yields the
two other inequalities.
The equivalence for {\bf (b)} holds.

For {\bf (c)}, we start with the remark that
$$
\{L_{a}\ge n\}=\{\exists k\ge n,\,|U_k|\ge a\}=\{\sup_{k\ge
n}|U_k|\ge a\}.
$$
The first equality is the definition of $L_{a}$, 
the second equality  is a consequence of the
almost sure convergence of $U_n$ to zero. (Hence, the
supremum is in fact almost surely a maximum.)
Write
$L_{a}'$,
resp.\ $L^*_{a}$, for
the functional $L_{a}$ associated to 
$(X_n')_n$,
resp.\ to
$(X^*_n)_n$, rather than to $(X_n)_n$.
We rely on the following inclusions:
\begin{eqnarray*}
& & \{L^*_{2a}\ge n\}\subset\{L_{a}\ge n\}
\cup\{L_{a}'\ge n\},
\\
\mbox{and}& & \{L_{2a}\ge n\}\cap\{L_{a}'\le
n\}\subset 
\{L^*_{a}\ge n\}.
\end{eqnarray*}
The first inclusion yields $P[L^*_{2a}\ge n]\le 2P[L_{a}\ge n]$, and
$$
E[G(L^*_{2a})]\le 2\,E[G(L_{a})].
$$
In the second inclusion,
$L_{2a}$ and $L_{a}'$ are independent and 
$L_{a}'$
is a.s.\ finite, hence $P[L_{a}'\le
n]$ converges to $1$ and $P[L_{2a}\ge n]$ is equivalent to
$P[L^*_{a}\ge n]$ when $n\to\infty$.
This proves that 
$G(L_{2a})$ is integrable as soon as $G(L^*_{a})$ is.
Finally, {\bf (c)} for $X^*$ holds if and only if
{\bf (c)} for $X$ does.
This ends the proof of proposition~\ref{th4}.

\section{The non moderate case}
\label{s:nondv}

The proof of theorem~\ref{th0} for non moderate functions is as follows.
If {\bf (a)} and {\bf (c)} are equivalent,
$G(L_{a})$ 
is integrable for all positive $a$ as soon as $|X|\,G(|X|)$ is integrable.
But $L_{2a}$ for $(2X_n)_n$ is $L_{a}$ for $(X_n)_n$,
hence
$|X|\,G(2|X|)$ should be integrable as soon as $|X|\,G(|X|)$ is.

Let $G$ be a non decreasing, increasing to infinity, non moderate function.
Then, $G(2t_n)\ge
n\,G(t_n)$ for an increasing sequence $(t_n)_n$.
Let the law of $X$ have as support the set $\{\pm t_n\,;\,n\ge 1\}$,
and give  to $t_n$ and $-t_n$ the same weight of
$c/(n^2\,t_n\,G(t_n))$, for a given constant $c$. The sum of the
weights converges, hence $c$ can be chosen so as to get a probability
measure. Then, $E[|X|\,G(|X|)]$ is finite and $E[|X|\,G(2|X|)]$ 
is infinite.
Thus, {\bf (a)} and {\bf (c)} are not equivalent for $X$.


The following properties of moderate functions are direct consequences
of the definition and we state them without proof.

Power and logarithmic functions are moderate,
exponential functions are not.
The function $G(t)=t^{b}\log(t)^{c}\log\log
(t)^{d}\ldots$
is moderate if the first non zero exponent is positive.
If $G(t)\sim a\,t^c$ when $t\to\infty$, with $c$ nonnegative, then
$G$ 
is moderate.
If $G$ is moderate,
then $G(t)\le c\,t^{c}$, for $c$ large enough and $t\ge 1$.
There exists non moderate $G$, such that $G(t)\le t^c$, for any
positive $c$.
There exists non moderate $G$, which are differentiable and such
that $G'(t)\le t^c$, for any positive $c$.
For any moderate $G$, there exists a smooth moderate
$H$ and $c\ge 1$, such that $c^{-1}\,H\le G\le c\,H$.
One can choose $H$ such that $H(t+1)/H(t)\le 1+a/(t+1)$, for
$a$ large enough.
Finally, there exists moderate $G$, such that the limit set of
$\log G(t)/\log(t)$ as $t\to\infty$ is the interval $[0,c]$, for any
positive $c$.


\bigskip
\noindent
LaPCS -- Laboratoire de Probabilit\'es, Combinatoire et Statistique
\\
Universit\'e Claude Bernard Lyon 1
\\
Domaine de Gerland
\\
50, avenue Tony-Garnier
\\
69366 Lyon Cedex 07 (France)

\medskip

\noindent
{\tt Didier.Piau@univ-lyon1.fr}
\\
{\tt http://lapcs.univ-lyon1.fr}


\end{document}